%% file: CI_10_25.tex
\documentclass{amsart}
\usepackage{amsmath, amsthm, latexsym, amssymb, amsfonts, epsfig, epsf}
\usepackage{pgf,pgfarrows,pgfnodes,pgfautomata,pgfheaps}

\newenvironment{pf}{\proof[\proofname]}{\endproof}
\theoremstyle{plain}
\newtheorem{Th}{Theorem}[section]

\newtheorem{Cor}[Th]{Corollary}
\newtheorem{Prop}[Th]{Proposition}
\newtheorem{Lemma}[Th]{Lemma}
\numberwithin{equation}{section}
\numberwithin{figure}{section}
\theoremstyle{definition}
\newtheorem{Rem}[Th]{Remark}
\newtheorem{Ex}[Th]{Example}
\newtheorem{Def}[Th]{Definition}

\newcommand{\cl}[1]{\mathcal{#1}}

\newcommand{\Z}{\mathbb Z}
\newcommand{\R}{\mathbb R}
\newcommand{\F}{\mathbb F}
\newcommand{\Pp}{\mathbb P}
\newcommand{\G}{\Gamma}

\newcommand{\D}{\Delta}

\newcommand{\cC}{\cl C}

\newcommand{\cL}{\cl L}

\newcommand{\cO}{\cl O}

\newcommand{\T}{{\mathbb T}}
\newcommand{\Jac}{{J_f^{\,\T}}}

\newcommand{\K}{{\mathbb K}}

\newcommand{\ov}{\overline}

\newcommand{\sig}{\sigma}
\newcommand{\Sig}{\Sigma}

\newcommand{\res}{\operatorname{res}}
\newcommand{\Res}{\operatorname{Res}}

\newcommand{\spn}{\operatorname{span}}

\newcommand{\rs}[1]{Section~\ref{S:#1}}
\newcommand{\rl}[1]{Lemma~\ref{L:#1}}
\newcommand{\rp}[1]{Proposition~\ref{P:#1}}
\newcommand{\rr}[1]{Remark~\ref{R:#1}}

\newcommand{\rc}[1]{Corollary~\ref{C:#1}}
\newcommand{\rt}[1] {Theorem~\ref{T:#1}}
\newcommand{\rd}[1]{Definition~\ref{D:#1}}
\newcommand{\rf}[1]{Figure~\ref{F:#1}}

\begin{document}


\title{Toric complete intersection codes}
\thanks{Partially supported by NSA Young Investigator Grant H98230-10-1-0163}
\author{Ivan Soprunov}
\address{Department of Mathematics, Cleveland State University, 2121 Euclid Ave, Cleveland OH, USA}
\email{i.soprunov@csuohio.edu}



\begin{abstract} In this paper we construct evaluation codes on zero-dimensional 
complete intersections in toric varieties and  give lower bounds for their minimum distance. This
generalizes the results of Gold--Little--Schenck and Ballico--Fontanari who considered
evaluation codes on complete intersections in the projective space. 
\end{abstract}

\maketitle


\section{Introduction}

This work is inspired by the results of 
Gold, Little and Schenck \cite{GLSch}, and Ballico and Fontanari \cite{BF}
on evaluation codes on complete intersections in the projective space. 
Examples of evaluation codes include Reed--Muller codes on points in affine and projective spaces and 
Goppa codes on points in algebraic curves. Here is a general definition. Let $X$ be an algebraic variety over a finite field $\F_q$
and let $S=\{p_1,\dots, p_N\}$ be a finite set of $\F_q$-rational points of $X$. Furthermore, let $\cL$ be a finite-dimensional
space of regular functions over $\F_q$ defined on an open subset of $X$ containing $S$. This defines an {\it evaluation map}
$$\text{ev}_{S}:\cL\to(\F_q)^{N},\quad f\mapsto (f(p_1),\dots,f(p_N)).$$ 
Its image is a linear code $\cC_{S,\cL}$ of block length $N$.
In the situation when $X$ is a projective toric variety, the set $S$ is the algebraic torus $(\F_q^*)^n$, and 
$\cL$ is the space of linear sections of a Cartier 
divisor on $X$ we obtain what is called a {\it toric code}.
In this case $\cL$ is spanned by monomials whose exponents are lattice points in a convex lattice polytope.
The minimum distance for toric codes was studied in \cite{Ha1, Jo, LiSche, LiSchw, Ru, SoSo1,SoSo2}.

Duursma, Renter\'ia, and Tapia-Recillas considered the situation when $X=\Pp^n$, the set $S$
is an arbitrary zero-dimensional complete intersection in $\Pp^n(\F_q)$, and $\cL=\cL_a$ is the space of homogeneous
polynomials of degree $a$. Their paper \cite{DRT}  is concerned
with computing the dimension of the corresponding evaluation codes $\cC_{S,\cL_a}$.
Later Gold, Little, and Schenck \cite{GLSch}  found a very nice application of the Cayley--Bacharach theorem 
that gave a lower bound for the minimum distance of $\cC_{S,\cL_a}$,
generalizing the 2-dimensional  result of Hansen \cite{Ha:2003}.  They showed 
that the minimum distance satisfies
$$d(\cC_{S,\cL_a})\geq s-a+2,$$
where  $s=\sum_{i=1}^n d_i-(n+1)$ and $d_1,\dots, d_n$ are the degrees of the polynomials defining~$S$.  
Ballico and Fontanari \cite{BF}  then gave a significantly better bound 
$$d(\cC_{S,\cL_a})\geq n(s-a)+2,$$
which holds 
for complete intersections $S$ satisfying a ``generality" condition: no $n+1$ points of  $S$ lie on a hyperplane in~$\Pp^n$. 

In this paper we combine the two situations: $X$ is a projective toric variety, $S$ is a zero-dimensional complete
intersection in $X$, and $\cL$ is a space of global sections of a Cartier divisor on $X$. The corresponding
evaluation code we call a {\it toric complete intersection code}.
We give two lower bounds for the minimum distance of such codes: for sets $S$ with
and without a ``generality'' condition. Our bounds generalize the ones in \cite{GLSch} and \cite{BF}.
Although we largely adopted methods from these papers, the difficulty is that no analog of the Cayley--Bacharach theorem
for toric varieties is currently known.
It turned out that the Toric Euler--Jacobi theorem (\rt{E-J}) on global residues (which can be thought of as a weak
toric analog of the Cayley--Bacharach theorem, see \rc{E-J1}) provides enough information for applications to evaluation codes.

In our exposition we decided to use not the language of toric geometry but rather the more explicit 
language of Laurent polynomial systems and Newton polytopes. The relationship between the two 
is discussed in \rs{toric}. \rs{prelim} gives the necessary preliminaries and states the Toric Euler--Jacobi
theorem and its immediate applications. \rs{bounds} contains the main results on the minimum distance
of toric complete intersection codes: \rt{noass} does not use any additional assumptions,
and \rt{ass} assumes a certain ``generality'' property of $S$.  In \rs{geom} we give geometric conditions on the Newton
polytopes of polynomials defining $S$ which guarantee that this property holds when the coefficients of the
polynomials are generic. The paper concludes with applications and concrete examples in \rs{ex} and remarks
about further work.

\section{Preliminaries}\label{S:prelim}

\subsection{Evaluation Codes}\label{S:eval}

In this section we will define evaluation codes we will be dealing with throughout the paper.
First let us introduce some standard definitions and notation from the theory of Newton polytopes.
Let $\K$ be a field and $\overline{\K}$ be its algebraic closure.
Consider a Laurent polynomial $f\in\K[t_1^{\pm 1},\dots,t_n^{\pm 1}]$. 
Its {\it Newton polytope}  $P(f)$ is the convex hull of the exponent vectors of the monomials appearing in $f$. Thus we can write
$$f=\sum_{a\in P(f)\cap\Z^n}c_at^a,\quad\text{where }\  t^a=t_1^{a_1}\cdots t_n^{a_n},\ \ c_a\in\K.$$
Given a face $Q$ of $P(f)$ the {\it restriction} $f^Q$ is the Laurent
polynomial
$$f^Q=\sum_{a\in Q\cap\Z^n}c_at^a.$$

Next we define evaluation codes slightly adapted to our situation (see also \cite{Ha2, Li, TV} for various
constructions of evaluation codes).
Choose a finite subset  $S=\{p_1,\dots,p_N\}$ of $(\K^*)^n$ and 
a finite-dimensional subspace $\cL$ of $\K[t_1^{\pm 1},\dots,t_n^{\pm 1}]$.
Define the {\it evaluation map}
$$\text{ev}_{S}:\cL\to \K^N,\quad f\mapsto (f(p_1),\dots,f(p_N)).$$
The image of $\text{ev}_{S}$ is a linear code, called the {\it evaluation code}, which we denote by $\cC_{S,\cL}$.

In the paper we will be dealing with evaluation codes $\cC_{S,\cL}$ where  $\cL$ is a
space of Laurent polynomials and $S$ is a zero-dimensional  complete intersection of $n$ hypersurfaces
in a toric variety.  We postpone the toric geometry definition of $S$ until \rs{toric}. Instead, we
formulate this in terms of the theory of Newton polytopes.  We describe $S$ as the solution set of  a Laurent polynomial system satisfying three assumptions below.
 
Fix a collection of $n$-dimensional convex lattice polytopes $P_1,\dots,P_n$ in $\R^n$
and let $P=P_1+\dots+P_n$ be their Minkowski sum.
Consider $n$ Laurent polynomials $f_1,\dots, f_n$ over $\K$ with Newton
polytopes $P_1,\dots,P_n$ such that the system $f_1=\dots=f_n=0$ satisfies the following.

\noindent {\bf Assumptions:}
\begin{enumerate}
\item the system is {\it non-degenerate} with respect to $P$,
i.e. for every proper face $Q\subset P$ the restricted system $f_1^{Q_1}=\dots=f_n^{Q_n}=0$ has no solutions
in $(\overline{\K}^*)^{n}$, where $Q=Q_1+\dots+Q_n$, for unique faces $Q_i\subset P_i$;
\item  at each $p\in S$ the collection $(f_1,\dots,f_n)$ forms a system of local parameters, i.e. the 1-forms $df_1,\dots,df_n$ are linearly independent at $p$;
\item the solution set $S\subset (\overline{\K}^*)^n$ of the system consists of $\K$-rational points.
\end{enumerate}

Before describing the space $\cL$ we need to set some notation. 
For any set $A\subset \R^n$ we use $A_\Z$ to denote the set of lattice points
in $A$, i.e. $A_\Z=A\cap\Z^n$.  Also, we let  $P^\circ$ denote the interior of 
the polytope $P=P_1+\dots+P_n$.
Now let $A$ be any subset of $P^\circ$. Define 
$$\cL(A)=\spn_\K\{t^a\ |\ a\in A_\Z\}\subset \K[t_1^{\pm 1},\dots,t_n^{\pm 1}].$$

\begin{Def}\label{D:toricCIcode} Let $S$ be the solution set of a system $f_1=\dots=f_n=0$  with $n$-dimensional
Newton polytopes $P_1,\dots, P_n$ satisfying (1)--(3) above. Let the set $A$ lie in the interior $P^\circ$
of $P=P_1+\dots+P_n$. The  evaluation code $\cC_{S,\cL(A)}$ is called a {\it toric complete intersection code}.
We will denote it  simply by $\cC_A$. Furthermore, $d(\cC_A)$ will denote the minimum distance (the minimum
weight) of $\cC_A$. 
\end{Def}

\begin{Rem}
Although the above definition makes sense for arbitrary subsets $A$ of $P^\circ$, we may just as well restrict 
ourselves to the case of convex polytopes $A$. Indeed, 
the construction of the code depend on $A_\Z$ rather than on $A$ itself. Moreover, the bounds  on the minimum distance 
of $\cC_A$ which we prove in \rs{bounds} will not change if one replaces $A$ with the convex hull of $A_\Z$, whereas the dimension
of $\cC_A$ may, of course, only increase. 
\end{Rem}

\subsection{The Toric Euler--Jacobi theorem}\label{S:E-J}

Here we discuss the toric analog of the Euler--Jacobi theorem (\rt{E-J}) and its consequences. 
This theorem was first discovered by  by Khovanskii (see \cite{uspehi}) over the field of complex numbers. 
In \cite[Sec. 14]{Kunz} the first part of the theorem is proved over an arbitrary algebraically closed field.
The second part of the theorem is proved over fields of positive characteristic by \cite{JA}
under the condition that the $P_i$ have the same normal fan, but is currently unknown
in general. Nevertheless, the proofs of our main results will only use the first part of \rt{E-J}, so we
do not make any additional assumptions on the polytopes (with the exception of \rt{generic2}).

\begin{Def} Let $f_1,\dots,f_n\in \K[t_1^{\pm 1},\dots,t_n^{\pm 1}]$ be Laurent polynomials. The
Laurent polynomial 
$$\Jac=\det\left(t_j\frac{\partial f_i}{\partial t_j}\right)$$ 
is called the {\it toric Jacobian} of $f_1,\dots, f_n$.
\end{Def}

It is easy to see that the Newton polytope $P(\Jac)$ of the toric Jacobian lies in $P=P_1+\dots+P_n$, where $P_i=P(f_i)$.  
Also, assumption (3) in \rs{eval} implies $\Jac(p)\neq 0$ for every $p\in S$. 

\begin{Th}[\cite{uspehi}]\label{T:E-J}
Let $S$ be the solution set of a system $f_1=\dots=f_n=0$  with $n$-dimensional
Newton polytopes $P_1,\dots, P_n$ satisfying (1)--(3) above.
Let $P=P_1+\dots+P_n$ be the Minkowski sum. Then
\begin{enumerate}
\item for any $h\in \cL(P^\circ)$ we have
$\sum_{p\in S}{h(p)}/{\Jac(p)}=0$;
\vspace{.2cm}
\item  for
any function $\phi: S\to \K$ with $\sum_{p\in S}\phi(p)=0$ there exists
 $h\in \cL(P^\circ)$ such that $\phi(p)=h(p)/\Jac(p)$ for every $p\in S$.
\end{enumerate}
\end{Th}

Here is an immediate corollary from the theorem.

\begin{Cor}\label{C:E-J1} Any $h\in \cL(P^\circ)$ which vanishes at $|S|-1$ points of $S$ must
vanish at all points of $S$. 
\end{Cor}

The next result, known as the Bernstein--Kushnirenko theorem, provides the size of the solution set $S$ for systems
$f_1=\dots=f_n=0$ with given Newton polytopes $P_1,\dots, P_n$.

\begin{Th} \label{T:Be}Let a Laurent polynomial system $f_1=\dots=f_n=0$  with 
Newton polytopes $P_1,\dots, P_n$ have isolated solution set $S$ in $(\overline{\K}^*)^n$.
Then $|S|$ cannot exceed the normalized mixed volume $V(P_1,\dots,P_n)$ of the Newton polytopes.
Moreover, $|S|=V(P_1,\dots,P_n)$ if and only if the system satisfies 
assumptions (1)-(2). 
\end{Th}

The original proof by \cite{Be} uses the homotopy continuation method and is valid over the field of complex numbers.
Kushnirenko in \cite{Kush}  gave an algebraic proof which works over any algebraically closed field  
regardless of the characteristic. A similar argument also appears in \cite[Sec. 6]{Tuit}.

\begin{Rem}\label{R:Be} Suppose we have a system $f_1=\dots=f_n=0$ with Newton polytopes $P_1,\dots, P_n$.
According to \rt{Be}, if we can exhibit $V(P_1,\dots,P_n)$-many $\K$-rational solutions to the system and the
solutions are isolated then the system must satisfy assumptions (1)--(3).
We will use this observation when constructing toric complete intersection codes in  \rs{ex}.
\end{Rem}

Here is our first application to toric complete intersection codes.

\begin{Prop}\label{P:E-J2} If $|S|>1$ then the minimum distance of $\cC_{P^\circ}$ is at least 2.
\end{Prop}

\begin{pf} Any  $h\in \cL(P^\circ)$ which is not identically zero on $S$ 
may have at most $|S|-2$ zeroes by \rc{E-J1}. Hence the weight of every non-zero codeword in $\cC_{P^\circ}$
is at least~2. To see that such $h$ exist one can show that if $|S|=V(P_1,\dots,P_n)>1$ then $P^\circ$ must
contain at least one lattice point $u$, and so $\cL(P^\circ)$ contains $t^u$. In fact, $V(P_1,\dots,P_n)=1$
is equivalent to all $P_i$ being equal to a basis simplex $\Delta$, in which case $P=n\Delta$ has no lattice points
\cite[Prop. 2.7]{CCDDS}.
\end{pf}

\subsection{Relation to Toric varieties}\label{S:toric}
Here we will show how our problem can be reformulated in the
language of toric geometry. Let $X=X_\Sig$ be a projective simplicial toric variety over $\K$ of dimension $n$,
defined by a complete rational simplicial fan $\Sig\subset\R^n$. Each ray $\rho\in \Sig(1)$ 
is generated by a primitive lattice vector $v_\rho\in\Z^n$ and corresponds to 
a torus-invariant prime divisor $D_\rho$ on~$X$.
A  {\it semi-ample} divisor $D$ on $X$ is a torus-invariant 
Cartier divisor $D=\sum_{\rho\in\Sig(1)}a_\rho D_\rho$
for which the corresponding line bundle $\cO(D)$ is generated by global sections. 
This implies that the set
$$P_D=\{u\in\R^n\ |\ \langle u,v_\rho \rangle\geq -a_\rho, \rho\in \Sig(1)\}$$
is a lattice polytope in $\R^n$ \cite[Sec. 3.4]{F}. Also the space of global $\K$-sections  of $\cO(D)$
is isomorphic to $\cL(P_D)$ in our notation in \rs{eval}.

Now fix $n$ semi-ample divisors $D_1,\dots, D_n$ on $X$ and let $P_i=P_{D_i}$ be the corresponding
lattice polytopes. Let $D=D_1+\dots+D_n$. 
For every $1\leq i\leq n$ let $f_i$ be a section of the line bundle $\cO(D_i)$.
The assumption (1) in \rs{eval} guarantees that the hypersurfaces  defined by the $f_i$ in $X$ do not
have common points on the orbits of $X$ of codimension greater than 1, which implies that the
hypersurfaces intersect in isolated points $S$ in the dense orbit. The other two assumptions say that the
intersections are transverse and consist of $\K$-rational points. 

The following is a higher-dimensional generalization of the $\Omega$-construction of evaluation
codes on algebraic curves \cite[Sec. 4.1.1]{TV}. Let $\Omega^n_X$ be the sheaf of Zariski  $n$-forms on $X$
and $\Omega^n_X(D)$ the sheaf corresponding to the divisor $D=D_1+\dots+D_n$. The global sections
of this sheaf are $n$-forms whose only  poles are in the support of the $D_i$. There is an isomorphism
$\Omega^n_X(D)\cong \cO(D-\sum_\rho D_\rho)$ \cite[Sec. 8.2]{CLSch}.  We can write this explicitly in 
affine coordinates $(t_1,\dots, t_n)$. A section of $\Omega^n_X(D)$ has the form
$$\omega_h=\frac{h}{f_1\cdots f_n}\frac{dt_1}{t_1}\wedge\cdots\wedge\frac{dt_n}{t_n},$$
for some Laurent polynomial $h$ which corresponds to a section of $\cO(D-\sum_\rho D_\rho)$.
Using the above identification, we see that the space of global sections of $\cO(D-\sum_\rho D_\rho)$  
is spanned by the lattice points of the (rational) polytope corresponding to $D-\sum_\rho D_\rho$, i.e.
the interior lattice points of $P_D=P_1+\dots+P_n$. Hence,  $h\in\cL(P^\circ)$.

Now let $S=\{p_1,\dots, p_N\}$ be the intersection of the hypersurfaces defined by the $f_i$ as above. Then
at every $p\in S$ the local (Grothendieck)
residue $\res_p(\omega_h)$ is defined \cite{GKh2}. Choose a subspace $\cL$ of global sections of $\Omega^n_X(D)$.
This results in the {\it residue map}
$$\res_S: \cL\to \K^{N}, \quad \omega_h\mapsto \left(\res_{p_1}(\omega_h),\dots, \res_{p_N}(\omega_h)\right),$$
whose image is a linear code.  In the case of transverse intersections 
at $p$ we have $\res_p(\omega_h)=h(p)/\Jac(p)$ and the reside map becomes:
$$\res_S: \cL\to \K^{N}, \quad \omega_h\mapsto \left(\frac{h(p_1)}{\Jac(p_1)},\dots, \frac{h(p_N)}{\Jac(p_N)}\right).$$
The linear code it defines is equivalent to the toric complete intersection code from \rd{toricCIcode}.
A similar construction of {\it toric residue codes} appears in \cite{JA} in relation to quantum stabilizer codes. 

The sum of the local residues over $p\in S$ is the {\it global
residue} $\Res_f(h)$ of $h$ with respect to $f=(f_1,\dots,f_n)$. In these terms the first statement of \rt{E-J} 
says that the global residue of any $h\in\cL(P^\circ)$ equals zero. The global residue is closely related to
the toric residue \cite{Coxres} and was  studied in \cite{CCD, CaD, global}.


\section{Bounds for the minimum distance}\label{S:bounds}

Recall that the evaluation code $\cC_A$ is constructed by choosing a subset $A$ of $P^\circ$.
Note that lattice translations of $A$, i.e. translations by lattice vectors, result in equivalent codes, so the minimum distance 
$d(\cC_A)$ is independent of such translations. Consider a  ``complementary'' set $B$, for which $A+B\subseteq P^\circ$.
 It turns out that $d(\cC_A)$ is related to properties of the space $\cL(B)$ as \rt{recip} below shows.
 The following definition from classical algebraic geometry will be used throughout the paper.
 
 \begin{Def} We say that a finite set of points $T\subset (\K^*)^n$ {\it imposes independent
 conditions} on a space of Laurent polynomials $\cL$ if the evaluation map $\text{ev}_{T}:\cL\to\K^{|T|}$
 is surjective.
\end{Def}

\begin{Th}\label{T:recip}
Let $S$ be the solution set of a system $f_1=\dots=f_n=0$ satisfying assumptions  (1)--(3) above.
 Let $A$ and $B$ be two subsets of $\R^n$ such that $A+B\subseteq P^\circ$. If any 
$T\subseteq S$ of size $m$ imposes independent conditions on the space 
$\cL(B)$ then $d(\cC_A)\geq m+1$.
\end{Th}

\begin{pf} We need to show that any $h\in\cL(A)$,  not identically zero on $S$, vanishes at no more than $|S|-m-1$ points
of $S$. Assume there exist $h\in\cL(A)$ and a subset $Z\subset S$  of size $|S|-m$ such that $h$ vanishes on $Z$,
but $h(p)\neq 0$ for some $p\in S$. By our assumption $S\setminus Z$ imposes independent conditions on $\cL(B)$,
so there exists $g\in\cL(B)$ such that $g$ vanishes at every point of $S\setminus (Z\cup \{p\})$, but not at $p$. 
Now the polynomial $hg$ belongs to $\cL(A+B)\subseteq\cL(P^\circ)$ and vanishes at every point of $S$ but not at $p$,
which contradicts \rc{E-J1}.
\end{pf}

\begin{Rem} Consider a special case: $X=\Pp^n$,  $f_1,\dots, f_n$ are homogeneous polynomials
of degrees $d_1,\dots, d_n$; and $\cL(A)$ and $\cL(B)$ are subspaces of homogeneous polynomials of 
degrees $a$ and $s-a$, respectively, where $s=\sum_{i=1}s^nd_i-(n+1)$. 
In this case \rt{recip} follows from the Cayley--Bacharach theorem \cite{EGH}
and serves as the main tool in the proofs of the results of \cite{GLSch} and \cite{BF}. 
We would like to point out that no toric analog of the Cayley--Bacharach theorem is currently known,
however, the Toric Euler--Jacobi theorem is sufficient for our application to toric complete intersection codes.
\end{Rem}

Our next goal is to understand what sets $B$ satisfy the condition of the above theorem for some value of $m$.
Here is our first example.

\begin{Lemma}\label{L:simple} Let $B=B_1+\dots+B_m$ where the lattice set  $B_i\cap\Z^n$
affinely generates $\Z^n$ for every $1\leq i\leq m$. Then any $m+1$ points in $(\K^*)^n$ impose
independent conditions on the space~$\cL(B)$.
\end{Lemma}

\begin{pf}   Suppose $m=1$ and let $T=\{p_0,p_1\}$ be any subset in $(\K^*)^n$.
It is enough to show that there is a polynomial  $g\in \cL(B)$ such that $g(p_1)=0$
and $g(p_0)\neq 0$.  We may assume that $B$ contains the origin.
Let $\{v_1,\dots,v_n\}\subseteq B$ be a basis for~$\Z^n$ and let $s=t^M=(t^{v_1},\dots,t^{v_n})$
be the corresponding automorphism of $(\K^*)^n$.  Choose a linear function $l(s)$ such that $l(p_1^M)=0$ and
$l(p_0^M)\neq 0$. Then  the polynomial $g(t)=l(t^M)$ lies in $\cL(B)$ and satisfies the required property.

In general,  let $T=\{p_0,\dots,p_m\}$ be any subset of $m+1$ points in $(\K^*)^n$.
By the previous case for every $1\leq i\leq m$ there exists
$g_i\in\cL(B_i)$ such that $g_i(p_i)=0$ and $g_i(p_0)\neq 0$. Then the polynomial 
$g=\prod_{i=1}^mg_i$ lies in $\cL(B)$, vanishes on $T\setminus\{p_0\}$,
and is not zero at~$p_0$. This implies that $T$ imposes independent conditions on $\cL(B)$.
\end{pf}

In our first application of \rt{recip} we estimate $d(\cC_A)$ using the number of ``primitive" simplices~$\D_i$ one can add to $A$ and
still stay in  $P^\circ$, after a possible lattice translation. We say that a simplex $\D$ is {\it primitive} if 
$\D=\text{conv.hull}\,\{0,v_1,\dots,v_n\}$, where $\{v_1,\dots,v_n\}$ is a basis for~$\Z^n$.

\begin{Th}\label{T:noass} Let $S$ be the solution set of a system $f_1=\dots=f_n=0$ satisfying assumptions  (1)--(3) above.
Let $A$ be any set such that $A+\D_1+\dots+\D_m\subseteq  P^\circ$ up to a lattice translation, where each $\D_i$ is a primitive simplex. Then $d(\cC_A)\geq m+2$.
\end{Th}

\begin{pf} This follows from \rt{recip} and \rl{simple}.
\end{pf}


In our next application we will consider solution sets $S\subset(\K^*)^n$ satisfying one additional assumption.

\smallskip

\noindent {\bf Assumption:}
\begin{enumerate}
\item[(4)] There exists an $n$-polytope $Q$ such that any $|Q_\Z|$
points of $S$ impose independent conditions on $\cL(Q)$. In other words,  for any subset $T\subset S$ of size  $|Q_\Z|$
the evaluation map  $\text{ev}_{T}:\cL(Q)\to\K^{|Q_\Z|}$ is an isomorphism.
\end{enumerate}

\begin{Ex} Suppose $X=\Pp^n$ and $Q=\D$ is the standard $n$-simplex, i.e. the convex hull of the origin and the
$n$ standard basis vectors. Then $(4)$ is equivalent to saying that no $n+1$ points of $S$ lie on a hyperplane.
Complete intersections in $\Pp^n$ with this ``generality" assumption were considered in \cite{BF}.
\end{Ex}

The assumption $(4)$ allows us to obtain better bounds on the minimum distance of the codes $\cC_A$, as was
suggested by \cite{BF} in the case of the projective space. In fact, their
approach generalizes to arbitrary toric varieties.
We will begin with a toric analog of their Horace Lemma.

\begin{Prop}\label{P:Horace} Let $T\subset(\K^*)^n$ be a finite subset and $A$ a bounded subset of $\R^n$.
Consider a hypersurface $H$ in $(\K^*)^n$ defined by $h\in\cL(Q)$.
If $T\cap H$ imposes independent conditions on $\cL(A+Q)$ and $T\setminus (T\cap H)$
imposes independent conditions on $\cL(A)$  then $T$ imposes independent conditions on $\cL(A+Q)$.
\end{Prop}

\begin{pf} Take any point $p\in T$. 
If $p\not\in H$ then there exists $g\in\cL(A)$ which does not vanish at $p$, but vanishes at all the other points 
of $T\setminus (T\cap H)$. Then the polynomial $f=gh\in\cL(A+Q)$ vanishes at all points of $T\setminus \{p\}$.
Also $f(p)=g(p)h(p)\neq 0$ since $p\not\in H$.

Now if  $p\in H$ then there exists $f_1\in\cL(A+Q)$ which does not vanish at $p$, but vanishes at all 
the other points of $T\cap H$. Consider 
the function $\phi:T\setminus (T\cap H)\to \K$ given by $q\mapsto f_1(q)/h(q)$. We know that there exists $g\in\cL(A)$
such that $g(q)=\phi(q)$ for any $q\in T\setminus (T\cap H)$. Put $f=f_1-gh$. Clearly $f\in\cL(A+Q)$ and
$f$ vanishes at every point of $T$ except at $p$.
\end{pf}

\begin{Prop}\label{P:induction} Let $S$ be any subset of $(\K^*)^n$ satisfying assumption (4).
Then, for any $k\geq 0$, any subset $T\subseteq S$ of size $|T|=(|Q_\Z|-1)k+1$  imposes independent conditions on $\cL({kQ})$.
\end{Prop}

\begin{pf} The proof is by induction on $k$. For $k=0$ we have $T=\{p\}$ which imposes independent
conditions on the space $\cL({kQ})\cong \K$.  

For $k>0$ choose  $T'\subset T$ of size $m=|Q_\Z|-1$.
Since $m<|Q_\Z|=\dim\cL(Q)$ there exists a non-zero polynomial $h\in\cL(Q)$ which vanishes on $T'$.
Moreover, $T'=S\cap H$, where $H$ is the hypersurface defined by $h$. Indeed, if 
$S\cap H$ contains a point $p$ not in $T'$ then the evaluation map  $\text{ev}_{T'\cup\{p\}}:\cL(Q)\to\K^{m+1}$ is 
degenerate which contradicts  the assumption~(4). Clearly, since $T'\subset T\subset S$ we have
$T'=S\cap H=T\cap H$.

Now $T\setminus T'$ has size $m(k-1)+1$ and by induction imposes independent conditions on  $\cL({(k-1)Q})$.
Also by $(4)$ the set $T'$ imposes independent conditions on $\cL(Q)$ and hence on $\cL({kQ})$ as $Q\subset kQ$ up to a lattice translation.
It remains to apply \rp{Horace}.
\end{pf}

\begin{Th} \label{T:ass}
Let $S$ be the solution set of a system $f_1=\dots=f_n=0$ satisfying assumptions (1)--(4).
Let  $A$ be any set such that $A+kQ\subset P^\circ$ up to a lattice translation, for some $k\geq 0$.
Then
$$d(\cC_A)\geq (|Q_\Z|-1)k+2.$$
\end{Th}

\begin{pf}  The theorem follows from  \rp{induction} and \rt{recip} where we put $m=(|Q_\Z|-1)k+1$.
\end{pf}

\section{Constructing toric complete intersection codes}\label{S:geom}

In this section we give geometric conditions on the polytopes  $P_1,\dots, P_n$ and $Q$ that produce systems 
 satisfying assumption $(4)$ if the coefficients are generic elements of $\overline{\K}$. We use these
 conditions when constructing examples of toric complete intersection codes in \rs{ex}.


\begin{Th}\label{T:generic} Let $Q$ be an $n$-dimensional lattice  polytope such that $Q_\Z$ generates $\Z^n$.
Suppose 
\begin{enumerate}
 \item[1.] $V(P_1,\dots,P_{n-1},Q)\geq |Q_\Z|$, 
 \item[2.] $(|Q_\Z|-1)Q\subset P_n$. 
 \end{enumerate}
 Then the solution set 
of any system $f_1=\dots=f_n=0$ with Newton polytopes $P_1,\dots, P_n$ and generic coefficients satisfies assumption $(4)$.
\end{Th}

\begin{pf}  Let $m=|Q_\Z|-1$. Let $\G_i$ be the hypersurface in $(\ov\K^*)^n$ defined by $f_i$.
Consider the curve  $C=\G_1\cap\dots\cap \G_{n-1}$ in $(\ov\K^*)^n$.  
Let  $V$ consist of all ordered collections  $(p_0,\dots, p_m)$ of regular points in $C$ such that $\{p_0,\dots, p_m\}$ do not impose independent conditions on $\cL(Q)$. In other words, 
$$V=\{T=(p_0,\dots, p_m)\in C_{\text{reg}}^{m+1}\ |\ \text{ev}_T:\cL(Q)\to(\ov\K^*)^{m+1} \text{ is not surjective}\,\},$$
where by abuse of notation we denote by $T$ both the ordered collection $(p_0,\dots, p_m)$ and the set $\{p_0,\dots, p_m\}$.
The set $V$ is algebraic with a dense open subset $V_0\subset V$ consisting of points of $V$  for which the map $\text{ev}_T$ has 
one-dimensional kernel.  

First we will show that $\dim V=m$. 
Indeed, every $T\in V_0$ defines  a unique hypersurface~$H$, defined by a polynomial in $\cL(Q)$,
such that the corresponding set $T$ lies in $C\cap H$.  We obtain a map $\pi:V_0\to\Pp\cL(Q)$.
On the other hand, by the Bernstein--Kushnirenko theorem (see \rt{Be}) 
any generic hypersurface $H$ with Newton polytope $Q$ satisfies
$|C\cap H|=V(P_1,\dots, P_{n-1},Q)\geq m+1$, so the image of $\pi$ is dense in $\Pp\cL(Q)$. Clearly, the fibers 
$\pi^{-1}(H)$ are finite, so we get $\dim(V)=\dim(V_0)=\dim(\pi(V_0))=\dim(\Pp\cL(Q))=m$. 

Now we will show that choosing a generic $f_n$ with Newton polytope $P_n$ produces 
$S=C\cap\G_n$ which satisfies assumption $(4)$. For this consider the set
$$W=\bigcup_{T\in V}W_T,\quad \text{where }\ W_T=\{f\in\cL(P_n)\ |\ f\text{ vanishes on } T \}.$$
Clearly, every $f_n$ in the complement of $W$  produces such $S$ (we also must avoid those $f_n$ 
which have zero coefficients corresponding to the vertices of $P_n$), so we need to show that $W$ has positive
codimension in $\cL(P_n)$. Indeed, according to our assumption $mQ\subset P_n$, so every set of $m+1$ points in $S$
imposes independent conditions on  $\cL(mQ)$ (by \rl{simple}) and hence on $\cL(P_n)$. Therefore the codimension of every subspace $W_T$ equals $m+1$. Thus $W$ is a vector bundle with $m$-dimensional base and codimension $m+1$ fibre,
so $W$ has codimension one.
\end{pf}


In the next theorem we show that in some situations the condition $(|Q_\Z|-1)Q\subset P_n$ can be replaced with 
$P_1+\dots+P_{n-1}+Q\subset P_n$. 
When $|Q_\Z|$ grows fast as a function of  $n$, the latter condition is preferable if one wants to 
avoid dealing with unnecessarily large $P_n$. 

We will need the following consequence of the Toric Euler--Jacobi theorem. 

\begin{Prop}\label{P:interpol}  Let  $P_1,\dots, P_n$
be  $n$-dimensional lattice polytopes  with the same normal fan, such that 
$\text{\rm char}\,\K$ does not divide the normalized mixed volume $V(P_1,\dots,P_n)$.
Let $S$  be the solution set for a system $f_1=\dots=f_n=0$  with Newton polytopes  $P_1,\dots, P_n$,
satisfying assumptions (1)--(3). Then $S$ imposes independent conditions on the space~$\cL(P)$.
\end{Prop}

\begin{pf} We need to show that for any function $\psi: S\to\K$ there exists $g\in\cL (P)$ with $g(p) =\psi (p) $
for all $p\in S$. Define $\phi:S\to\K$ by setting $\phi (p )= \frac{\psi(p )}{\Jac(p )}-c$, where 
$c=\frac{1}{|S|}\sum_{p \in S}\frac{\psi(p )}{\Jac(p )}$. Then $\sum_{p\in S}\phi(p )=0$, so by \rt{E-J}
there exists $h\in\cL(P^\circ)$ such that $h(p )=\Jac(p )\phi(p )$ for all $p\in S$. Now we can put $g=h+c\Jac\in\cL(P )$,
as $g(p )=h(p )+c\Jac(p )=\psi(p )$ for all $p\in S$, as required.
\end{pf}

This can be slightly refined. As we have seen in the above proof, \rp{interpol} still holds if
we replace $\cL(P)$ with $\spn_{\K}\{\cL(P^\circ),\Jac\}$.  

\begin{Th}\label{T:generic2} Let $P_1,\dots, P_n$ and $Q$ be $n$-dimensional lattice  polytopes with the same
normal fan and  such that $Q_\Z$ generates $\Z^n$. Suppose 
\begin{enumerate}
 \item[1.] $V(P_1,\dots,P_{n-1},Q)\geq |Q_\Z|$, 
 \item[2.] $P_1+\dots+P_{n-1}+Q\subset P_n$. 
 \end{enumerate}
 Then the solution set 
of any system $f_1=\dots=f_n=0$ with Newton polytopes $P_1,\dots, P_n$ and generic coefficients satisfies assumption $(4)$.
\end{Th}

\begin{pf} The proof is the same as for \rt{generic}, except for the last two sentences. Instead we need the following
observation.  Let  $T\in V$. By the definition of $V$ there exists 
a hypersurface $H$ defined by a polynomial in $\cL(Q)$ such that $T\subseteq C\cap H$. By
\rp{interpol},  $C\cap H$  imposes independent conditions on the space 
$\cL(P_1+\dots+P_{n-1}+Q)$. Therefore $T$ imposes independent conditions on $\cL(P_n)$ and hence
the codimension of the subspace $W_T$ equals $m+1$. The rest is as in the proof of \rt{generic}.
\end{pf}

\section{Examples}\label{S:ex}

In this section we put several applications of the results of the previous section as well as provide specific examples
of toric complete intersection codes over finite fields.

We start by showing how  \rt{noass} and \rt{ass}  recover the results of Gold--Little--Schenck and Ballico--Fontanari \cite{GLSch, BF}.

\begin{Ex} Let $S$ be a zero-dimensional smooth complete intersection in $\Pp^n$ given by
 $n$ homogeneous polynomials $F_1,\dots, F_n$ over $\K$.  Suppose $S$ lies in $\Pp^n(\K)$.
 Up to a projective change of coordinates we may assume that $S$ lies in the algebraic torus $(\K^*)^n$.
Rewriting  $F_i$ in the affine coordinates for $(\K^*)^n$ we obtain a polynomial $f_i$
with Newton polytope $P_i=d_i\D$ where $\D$ is the standard $n$-simplex and $d_i=\deg(F_i)$.
It is easy to see that $S$ satisfies the assumptions (1)--(3) in \rs{eval}. 

Now let $s=\sum_{i=1}^nd_i-(n+1)$ and let $A=a\D$ for some $1\leq a\leq s$. 
Notice that $\cL(A)$ is the space of polynomials of total degree at most $a$. We are going to apply
\rt{noass} with $m=s-a$ and all the $\D_i$ being simply $\D$. 
Clearly,  $A+\D_1+\dots+\D_n$, which equals $s\D$, lies in the interior of  $P=(\sum_{i=1}^nd_i)\D$.
Therefore, by \rt{noass}, $d(\cC_A)\geq s-a+2$. This is the result of \cite{GLSch}.

Next suppose $S$ satisfies assumption (4) with $Q=\D$. As pointed out before this means that no $n+1$ points
of $S$ lie in a hyperplane in $\Pp^n$. Applying \rt{ass} with $k=s-a$ we obtain $d(\cC_A)\geq n(s-a)+2$,
which is the result of \cite{BF}.
\end{Ex}

\medskip

In the next example we consider systems defined by multi-homogeneous polynomials.
This is the case of toric variety $X=\Pp^1\times\dots\times\Pp^1$.

\begin{Ex} For $1\leq i\leq n$ let $P_i$ be the lattice box with dimensions $(d_{i1},\dots,d_{in})$, each $d_{ij}\geq 1$. 
Let $S$ be the solution set of a system $f_1=\dots=f_n=0$ with Newton polytopes $P_1,\dots, P_n$ satisfying
assumptions (1)--(3).  
By the Bernstein--Kushnirenko theorem $|S|=V(P_1,\dots,P_n)$ which equals $\text{Perm}(D)$, the permanent of
the matrix $D=(d_{ij})$. Indeed, since each $P_i$ is the Minkowski sum of segments $P_i=\sum_{j=1}^n I_{ij}$, where
$I_{ij}=[0,d_{ij}e_j]$, by the multi-linearity of the mixed volume we obtain
\begin{eqnarray}
V(P_1,\dots,P_n)=V\Big(\sum_{j=1}^nI_{1j},\dots,\sum_{j=1}^n I_{nj}\Big)
&=&\sum_{\sig\in \text{S}_n}V\big(I_{1\sig(1)},\dots,I_{n\sig(n)}\big)\nonumber\\
&=&\sum_{\sig\in \text{S}_n}d_{1\sig(1)}\cdots d_{n\sig(n)}=\text{Perm}(D)\nonumber.
\end{eqnarray}

Now let $A$ be a lattice box with dimensions $(a_1,\dots, a_n)$. Note that $P$ is a lattice box with
dimensions $(d_1,\dots, d_n)$, where $d_j=\sum_i d_{ij}$. Hence $A$ lies in $P^\circ$ whenever $1\leq a_j\leq d_j-2$. 
Next, suppose $S$ satisfies the assumption (4) with $Q=\square$, the unit $n$-cube. Then for $k=\min_{j}(d_j-2-a_j)$ we have
$A+k\square\subset P^\circ$. Applying \rt{ass} we get 
$$d(\cC_A)\geq (2^n-1)\min_{1\leq j\leq n}(d_j-2-a_j)+2.$$

Let us now see under which condition on the polytopes $P_i$ the assumption (4) is generically satisfied. 
According to \rt{generic} and \rt{generic2} it is enough to require $V(P_1,\dots,P_{n-1},\square)\geq 2^n$ and 
either $(2^n-1)\square\subseteq P_n$ or $P_1+\dots+P_{n-1}+\square\subseteq P_n$.
The latter occurs when $d_{nj}\geq \min(2^n-1,\sum_{i=1}^{n-1}d_{ij}+1)$ for $1\leq j\leq n$. 
For the former note that $\square\subset P_i$, so by monotonicity of the mixed volume
$$V(P_1,\dots,P_{n-1},\square)\geq V(\square,\dots,\square)=n!\geq 2^n,$$
for $n\geq 4$. For $n=2$ we require $V(P_1,\square)=d_{11}+d_{12}\geq 4$. For $n=3$ we require that
at least one edge of either $P_1$ or $P_2$ has length 2, since in this case
$$V(P_1,P_2,\square)=d_{11}d_{22}+d_{12}d_{23}+d_{13}d_{21}+d_{13}d_{22}+d_{12}d_{21}+d_{11}d_{23}\geq 8.$$
\end{Ex}

\medskip

In the next two examples we present two explicit toric complete intersection codes over  $\F_{16}$ and $\F_{128}$,
respectively. In both cases the toric variety is a del Pezzo surface. We use MAGMA \cite{magma} for
constructing these examples. 


\begin{Ex} Let $\xi$ be a generator of the cyclic group $\F_{16}^*$. Let $P_1$ and $P_2$ be as in \rf{ex3}.
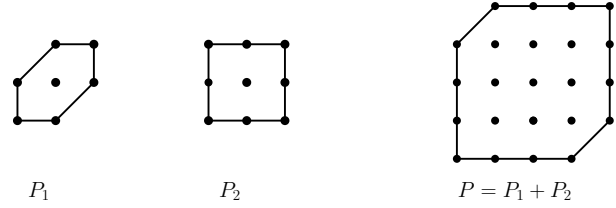
\begin{figure}[h]
\centerline{
 \scalebox{0.4}
 {
\input{ex3.pstex_t}}}
\caption{The Newton polygons and their Minkowski sum}
\label{F:ex3}
\end{figure}
Consider the following system. 
$$\begin{cases} f_1=x^2y^2 +  \xi^{14}x^2y +  \xi^{11}xy^2 +  \xi^9xy +  \xi^4x +  \xi^7y +  \xi^2=0,\\ 
                              f_2=x^2y^2 + \xi^8x^2y + \xi^4x^2 + \xi^4xy^2 + \xi xy + \xi^7x + \xi^{11}y^2 + \xi^7=0.\end{cases}$$
The system has $8=V(P_1,P_2)$ simple solutions in $(\F_{16}^*)^2$:
$$S=\{(\xi^7, \xi^6), (\xi^8, \xi), (\xi^{9}, 1), (\xi^{10}, \xi^{9}), (\xi^{12}, \xi^{14}), (\xi^{13}, \xi^8), (\xi^{14}, \xi), (\xi^{14}, \xi^{11})\}.$$
Let $Q=\square$, the unit square. 
One can check that any 4 points of $S$ impose independent conditions on the space $\cL(\square)$.
Now choose $A=\square$ as well. We have $A+\square\subset P^\circ$, so $$d({\cC}_{A})\geq (4-1)+2=5.$$
Furthermore $\dim{\cC}_{A}=|A_{\Z}|=4$, so we get an MDS $[8,4,5]$-code over $\F_{16}$.
\end{Ex}

To construct a bigger example we start with polygons $P_1$, $P_2$ satisfying the conditions of \rt{generic}. Then
we choose a random polynomial $f_1$ with Newton polytope $P_1$.  If the size of the field is big enough we can 
choose $V(P_1,P_2)$ rational points on the curve $f_1=0$ which satisfy assumption (4).

\begin{Ex}\label{Ex:4} 
The polygons $P_1$ and $P_2$ and their Minkowski sum $P$ are depicted in \rf{ex4}.  Consider
a system $f_1=f_2=0$ over $\F_q$ with Newton polytopes $P_1,P_2$ satisfying assumptions (1)--(3),
and let $S$ be the solution set of the system. 
We have $|S|=V(P_1,P_2)=14$. On the other hand,  a simple application of the Serre bound
shows that  for $q\leq 8$ the curve $f_1=0$ has less than 14 rational points. Therefore we must have $q>8$.

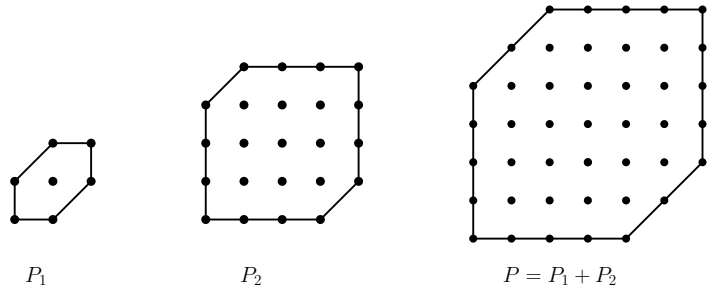
\begin{figure}[h]
\centerline{
 \scalebox{0.4}
 {
\input{ex5.pstex_t}}}
\caption{The Newton polygons and their Minkowski sum}
\label{F:ex4}
\end{figure}

First we consider an application of \rt{noass}. 
Take $A$ to be a $2\times 2$ lattice square, $\D_1$ the convex hull of $\{(0,0), (1,0), (1,1)\}$,
and $\D_2$ the convex hull of $\{(0,0), (0,1), (1,1)\}$. Then $A+\D_1+\D_2\subset P^\circ$. Therefore,
by \rt{noass}, we have $d(\cC_A)\geq 2+2=4$. The evaluation map 
$\text{ev}_{S}:\cL(A)\to{\F_{q}}^{\!14}$ has one-dimensional kernel spanned by $f_1$.
Therefore,  $\dim\cC_A=|A_{\Z}|-1=8$ and we obtain a $[14, 8, \geq 4]$-code over $\F_q$ with $q\geq 9$.

Next we consider a set $S$ satisfying the additional assumption (4). We set $Q=\square$, the unit square.
Since $V(P_1,\square)=4$ and $3\,\square\subset P_2$, both conditions of \rt{generic} are satisfied. 
 
We will work over $\F_{128}$; as before $\xi$ will denote a generator of $\F_{128}^*$.  
Here is a random polynomial over $\F_{128}$ with Newton polytope $P_1$:
$$f_1=x^2y^2 + \xi x^2y + \xi^{32}xy^2 + \xi^{4}xy + \xi^{78}x + \xi^{110}y + \xi^{31}.$$
The curve $f_1=0$ has 146 rational points in the torus. We choose 14 of these rational
points which impose independent conditions on $\cL(\square)$. Here is one such subset:
\begin{eqnarray}
&S&=\{(\xi^{5}, \xi^{91}), (\xi^{43}, \xi^{59}), (\xi^{44}, \xi^{100}), (\xi^{47}, \xi^{125}), (\xi^{51}, \xi^{33}), (\xi^{58}, \xi^{90}), (\xi^{68}, \xi^{42}), 
\nonumber\\
& & \ \ \ \ \, (\xi^{78}, \xi^{11}), (\xi^{78}, \xi^{12}), (\xi^{85}, \xi^{79}), (\xi^{96}, \xi^{11}), (\xi^{105}, \xi^{41}), (\xi^{116}, \xi^{106}), (\xi^{124}, \xi^{65})\}.
\nonumber
\end{eqnarray}
Since $|P_2\cap\Z^2|>|S|=14$ there exist polynomials $f_2$ with Newton polytope $P_2$ which vanish at $S$. 
We choose such $f_2$ that has no common factors with $f_1$. For example, we can take
\begin{eqnarray}
&f_2&=x^4y^4 +  \xi^{59}x^4y +\xi^{10}x^3y + \xi^{66}x^3 + \xi^{26}x^2y + \xi^{104}x^2 +  xy^4 + \xi^{44}xy^3 + \xi^{50}xy^2
\nonumber\\ 
& &\  \ \  \ \ \ \ \ \ \  + \xi^{78}xy + \xi^{56}x + \xi^{118}y^3 + \xi^{38}y^2 + \xi^{36}y + \xi^{108}.\nonumber
\end{eqnarray}
By \rr{Be}, $S$ is the solution set of $f_1=f_2=0$ and satisfies assumptions (1)--(4).
Next we look at different choices of the set $A$.

\begin{enumerate}
\item[(a)]  Let $A=P_1$. Then $A+2\,\square\subset P^\circ$, so
by \rt{ass} we get $d(\cC_A)\geq (4-1)\cdot 2+2=8$. It is easy to see that $\dim\cC_A=|A_{\Z}|-1=6$ and we obtain
a $[14,6,\geq 8]$-code over $\F_{128}$. In fact, the minimum distance is exactly 8 in this case.
\item[(b)] Let $A$ be the segment joining $(0,0)$ and $(1,1)$. Then $A+3\,\square\subset P^\circ$, so 
by \rt{ass} we get $d(\cC_A)\geq (4-1)\cdot 3+2=11$. Since $\dim\cC_A=|A_{\Z}|=2$ we get a
$[14,2,\geq 11]$-code over $\F_{128}$, which is in fact a $[14,2,13]$-code.
\item[(c)] Let $A=P_1+\square$.
Then $A+\square\subset P^\circ$ so
by \rt{ass} we get $d(\cC_A)\geq (4-1)+2=5$. To compute the dimension of $\cC_A$ note that $\text{ev}_S:\cL(A)\to{\F_{128}}^{\!14}$
has 4-dimensional kernel. In fact, $\cL(A)\cap I=\spn\{f_1, xf_1, yf_1, xyf_1\}$, where $I$ is the ideal
generated by $f_1, f_2$.  Therefore $\dim\cC_A=|A_\Z|-4=10$.  This shows that $\cC_A$ is an MDS code over $\F_{128}$ 
with parameters $[14,10,5]$.
\end{enumerate}

\end{Ex}

\section{Conclusion and further work}

Given a system of Laurent polynomial equations $f_1=\dots=f_n=0$ with $n$-dimensional Newton 
polytopes $P_1,\dots,P_n$ satisfying assumptions (1)--(3) or (1)--(4) and a set $A\subset P^\circ$
we  defined a  class of evaluation codes $\cC_{S,\cL(A)}$, called toric complete intersection codes, and
found general lower bounds for their minimum distance. 

We then gave conditions on the polytopes $P_1,\dots,P_n$ and $Q$ which guarantee that 
generic systems with such Newton polytopes satisfy assumption (4). One would like
to obtain some general results about the size of the field for which toric complete intersections
with given polytopes exist and with what probability they occur.  This would allow a more
systematic way of constructing them and studying their parameters.



Computing the dimension of $\cC_{S,\cL(A)}$ is not obvious since the evaluation map will have a non-trivial kernel, in general.
It requires computing the codimension of the ideal generated by the $f_i$ in the space $\cL(A)$. 
Although this can be done in concrete examples  one would like to have a general way of doing so.

\section*{Acknowlegements}
I thank  \c Stefan Tohaneanu for several fruitful discussions about evaluation codes on complete intersections
and explaining his work in \cite{To1, To2}; 
and Jan Tuitman for answering questions about the Bernstein--Kushnirenko theorem in positive characteristic.
I am grateful to two anonymous referees whose comments helped to improve the exposition.

\end{document}

%% file: ex3.pstex_t
\begin{picture}(0,0)%
\includegraphics{ex3.pstex}%
\end{picture}%
\setlength{\unitlength}{3947sp}%
\begingroup\makeatletter\ifx\SetFigFont\undefined%
\gdef\SetFigFont#1#2#3#4#5{%
  \reset@font\fontsize{#1}{#2pt}%
  \fontfamily{#3}\fontseries{#4}\fontshape{#5}%
  \selectfont}%
\fi\endgroup%
\begin{picture}(9449,3202)(523,-3193)
\put(7501,-3061){\makebox(0,0)[lb]{\smash{{\SetFigFont{20}{24.0}{\rmdefault}{\mddefault}{\updefault}{\color[rgb]{0,0,0}$P=P_1+P_2$}%
}}}}
\put(3751,-3061){\makebox(0,0)[lb]{\smash{{\SetFigFont{20}{24.0}{\rmdefault}{\mddefault}{\updefault}{\color[rgb]{0,0,0}$P_2$}%
}}}}
\put(751,-3061){\makebox(0,0)[lb]{\smash{{\SetFigFont{20}{24.0}{\rmdefault}{\mddefault}{\updefault}{\color[rgb]{0,0,0}$P_1$}%
}}}}
\end{picture}%

%% file: ex5.pstex_t
\begin{picture}(0,0)%
\includegraphics{ex5.pstex}%
\end{picture}%
\setlength{\unitlength}{3947sp}%
\begingroup\makeatletter\ifx\SetFigFont\undefined%
\gdef\SetFigFont#1#2#3#4#5{%
  \reset@font\fontsize{#1}{#2pt}%
  \fontfamily{#3}\fontseries{#4}\fontshape{#5}%
  \selectfont}%
\fi\endgroup%
\begin{picture}(10949,4477)(823,-3868)
\put(8551,-3736){\makebox(0,0)[lb]{\smash{{\SetFigFont{20}{24.0}{\rmdefault}{\mddefault}{\updefault}{\color[rgb]{0,0,0}$P=P_1+P_2$}%
}}}}
\put(1051,-3736){\makebox(0,0)[lb]{\smash{{\SetFigFont{20}{24.0}{\rmdefault}{\mddefault}{\updefault}{\color[rgb]{0,0,0}$P_1$}%
}}}}
\put(4426,-3736){\makebox(0,0)[lb]{\smash{{\SetFigFont{20}{24.0}{\rmdefault}{\mddefault}{\updefault}{\color[rgb]{0,0,0}$P_2$}%
}}}}
\end{picture}%